\documentclass[11pt,oneside]{article}
\usepackage{amssymb,amsmath}
\usepackage{url}
\renewcommand{\le}{\leqslant}
\renewcommand{\ge}{\geqslant}

\newcommand{\Span}{\mathrm{span}}

\newcommand{\RR}{\mathbb{R}}
\newcommand{\ZZ}{\mathbb{Z}}
\newcommand{\QQ}{\mathbb{Q}}
\newcommand{\MM}{\mathbf{M}}
\newcommand{\NN}{\mathbb{N}}

\newcommand{\CC}{\mathbb{C}}

\newcommand{\VVV}{\mathcal{V}}

\newcommand{\CCC}{\mathcal{C}}

\newcommand{\FFF}{\mathcal{F}}
\newcommand{\SSS}{\mathcal{S}}

\newcommand{\PPP}{\mathcal{P}}

\newcommand{\MMM}{\mathcal{M}}

\newcommand{\III}{\mathcal{I}}

\newcommand{\DDD}{\mathcal{D}}
\newcommand{\AAA}{\mathcal{A}}

\newcommand{\vx}{\mathbf{x}}

\newcommand{\vp}{\mathbf{p}}
\newcommand{\vq}{\mathbf{q}}

\newcommand{\va}{\mathbf{a}}
\newcommand{\vf}{\mathbf{f}}

\newcommand{\valp}{\boldsymbol{\alpha}}

\newtheorem{lemma}{Lemma}
\newtheorem{theorem}{Theorem}
\newtheorem{proposition}{Proposition}

\newtheorem{corollary}{Corollary}
\newtheorem{theoremb}{Theorem B}

\newtheorem{conjectureby}{Conjecture BY}

\newenvironment{proof}{\textsc{Proof. }}{\ \newline\hspace*{\fill}$\boxtimes$}

\newcommand{\vol}{\mathrm{Vol}}
\newcommand{\codim}{\mathrm{codim}\;}

\begin{document}
\title{Simultaneous Diophantine approximation on the three dimensional Veronese curve}

\author{Dmitry Badziahin}

%\thanks{The University of Sydney, dzmitry.badiahin@sydney.edu.au}

\maketitle

\begin{abstract}
We compute the Hausdorff dimension of the set of simultaneously
$\lambda$-well approximable points on the Veronese curve in $\RR^3$
for $1/3\le \lambda\le 3/5$. This range for $\lambda$ was predicted
in the conjecture of Beresnevich and Yang from~\cite{ber_yan_2023}.
To the best of the author's knowledge, this makes $\VVV_3$ the first
nondegenerate curve in $\RR^n$, $n\ge 3$, to confirm the lower bound
part of this conjecture.
\end{abstract}

{\footnotesize{{\em Keywords}: simultaneously well approximable
points, Veronese curve, Hausdorff dimension, simultaneous
Diophantine approximation on manifolds, cubic polynomials of bounded
discriminant

Math Subject Classification 2020: 11J13, 11J54, 11J82, 11K55}}

\section{Introduction}

For a positive real number $\lambda$ the set $S_n(\lambda)$ of simultaneously
$q^{-\lambda}$-well approximable points in $\RR^n$ is defined as follows:
$$
S_n(\lambda):= \{\vx\in\RR^n: ||q\vx - \vp||_\infty < q^{-\lambda}\;\mbox{ for i.m. } (q,\vp)\in \ZZ^{n+1}\}.
$$
One of the central problems in the metric theory of Diophantine
approximation is to understand the structure of the intersection of
$S_n(\lambda)$ with a suitable manifold $\MMM$. For more
information, the interested reader is referred to the influential
papers on the topic:~\cite{kle_mar_1998, beresnevich_2012,
ber_yan_2023}. We are particularly interested in computing the
Hausdorff dimension of this intersection. It is well known that this
dimension depends essentially on the choice of $\MMM$. For example,
consider the circle $\CCC:= \{\vx\in\RR^2: x_1^2+x_2^2=3\}$. It is
not difficult to verify that for $\lambda>1$ and $\vx\in \CCC$, the
inequality $||q\vx -\vp||_\infty < q^{-\lambda}$,
$(q,\vp)\in\ZZ^3\setminus \mathbf{0}$ for large $q$, implies that
$\vp/q\in\CCC$. From this, we immediately deduce that
$S_2(\lambda)\cap \CCC = \emptyset$. On the other hand, for the
Veronese curve $\VVV_n:=\{(x,x^2,\ldots, x^n): x\in\RR\}$ and all
$\lambda>1$ Schleischitz~\cite{schleischitz_2016} showed that $\dim
(S_n(\lambda)\cap \VVV_n) = \frac{2}{n(1+\lambda)}$, in sharp
contrast to the previous example.

However, if $\lambda$ is sufficiently close to $1/n$ and $\MMM$
satisfies certain natural nondegeneracy conditions then it is known
that $\dim(S_n(\lambda)\cap \MMM)$ does not depend on the choice of
the manifold. In particular, Beresnevich~\cite{beresnevich_2012}
showed that for $\lambda$ sufficiently close to $1/n$ and
nondegenerate $\MMM$,
$$
\dim (S_n(\lambda)\cap\MMM) \ge \dim \frac{n+1}{\lambda+1} - \codim\MMM.
$$
If $\MMM$ is one dimensional, that is, a nondegenerate curve, the
above inequality is achieved for $\frac{1}{n}\le\lambda\le
\frac{3}{2n-1}$. Later, in~\cite{ber_yan_2023, sst_2024} it was
shown that this lower bound is sharp when $\lambda$ is very close to
$1/n$, much closer than $\frac{3}{2n-1}$. While the general bound on
$\lambda$ in these results is rather intricate and depends on the
degree of nondegeneracy of $\MMM$, we state it here only in the case
of curves. Suppose that $\CCC$ is a curve parametrised by $n$ times
continuously differentiable function $\vf: J\to \RR^n$ such that its
derivatives up to degree $n$ span $\RR^n$ at any $x\in J$.
Then~\cite{sst_2024}
\begin{equation}\label{sst}
\dim (S_n(\lambda)\cap\CCC) = \frac{2-(n-1)\lambda}{1+\lambda}\qquad
\forall\; \lambda \in \left[\frac1n,
\frac{1}{n}+\frac{n+1}{n(2n-1)(n^2+n+1)}\right).
\end{equation}
Notice that the upper bound on $\lambda$ in~\eqref{sst} is of the form $\frac{1}{n}+O\big(\frac{1}{2n^3}\big)$.

A natural and interesting question then arises --- first formally
posed by Beresnevich and Yang~\cite{ber_yan_2023}. Let $\MM$ be a
class of $d$-dimensional manifolds in $\RR^n$. Define $\tau(\MM)$ to
be the supremum of values $\tau$ such that
$$
\dim(S_n(\lambda)\cap \MMM) = \frac{n+1}{\lambda+1} - \codim \MMM\quad\mbox{whenever}\; \frac1n\le \lambda< \tau
$$
for all manifolds $\MMM\in\MM$. The problem is then to compute, or
at least obtain nontrivial bounds for, $\tau(\MM)$ for a suitably
chosen class of manifolds $\MM$. The authors~\cite[Conjecture
2.7]{ber_yan_2023} proposed the following conjecture:
\begin{conjectureby}
Let $\MM_{n,1}$ be the set of nondegenerate curves in $\RR^n$; that
is, curves $\CCC=\{\vf(x): x\in I, \vf:I\to \RR^n\}$, such that they
are continuously differential enough times and that the set of their
derivatives at each point $x\in I$ spans $\RR^n$. Then
$$
\tau_{n,1} = \tau(\MM_{n,1}) = \frac{3}{2n-1}.
$$
\end{conjectureby}

In this paper, we continue the investigation initiated
in~\cite{badziahin_2025}. There, the focus was on the Veronese
curve: $\MMM = \VVV_n:=\{(x,x^2,, \ldots, x^n): x\in I\}$ where $I$
is an interval in $\RR$. For convenience, we use the notation
$S_n(I,\lambda):= S_n(\lambda)\cap \VVV_n$. In this setting, we are
able to substantially strengthen the known bounds for $\lambda$:

\begin{theoremb}
For all $\lambda$ between $\frac1n$ and $\frac{2}{2n-1}$ one has
$$
\dim S_n(I,\lambda) = \frac{2-(n-1)\lambda}{1+\lambda}.
$$
For $n=3$ the range for $\lambda$ can be extended to $\frac13\le
\lambda \le \frac{1}{2}$.
\end{theoremb}

Notice that the upper bound for $\lambda$ is asymptotically $\frac1n
+ O\big(\frac{1}{2n^2}\big)$ which is higher than that in the latest
result~\eqref{sst} from~\cite{sst_2024} for general nondegenerate
curves.

We focus on the case of $n=3$ and extend the range of $\lambda$ in
Theorem~B to $\lambda\le \frac35$. This upper bound coincides with
the value $\tau_{3,1}$ predicted in Conjecture BV. Therefore, we
confirm that the set of $\lambda$-well approximable points on
$\VVV_3$ behaves as predicted by the conjecture.

\begin{theorem}\label{th1}
For all $\lambda$ between $\frac13$ and $\frac35$ one has
$$
\dim S_3(I,\lambda) = \frac{2-2\lambda}{1+\lambda}.
$$
\end{theorem}

Throughout the paper, we use Vinogradov notation. For positive real
quantities $A$ and $B$, we write $A\ll B$ if $A\le cB$ for some
constant $c>0$ that may depend only on the manifold $\MMM$ (in
particular, on the dimension $n$ of the ambient space) but not on
the specific rational points $\vq$. The notion $A\gg B$ is defined
analogously and $A\asymp B$ means that both $A\ll B$ and $A\gg B$
hold simultaneously.

\section{General setup}

We begin by presenting several general techniques for estimating
upper bounds on the Hausdorff dimension of limsup sets. These
methods are not new and appear, often implicitly, in many existing
papers. Our main purpose in presenting them here is to provide a
clear point of reference. In the following sections, we will make
extensive use of these techniques. We also hope they will prove
useful in future work.

Let $(B_i)_{i\in \III}$ be a sequence of balls in a metric space $X$
equipped with a height function $q:\III\to \NN$. Suppose that for
each $q\in \NN$ the set $I(q):=\{i\in \III: q(i)=q\}$ is finite and
that there exists a function $\rho:\NN\to\RR^+$ such that
$\lim_{q\to\infty} \rho(q)=0$ and
$$
\forall i\in \III(q),\quad |B_i|\le \rho(q).
$$
Here, $|B|$ denotes the diameter of the set $B\subset X$.

We consider the problem of computing an upper bound for the
Hausdorff dimension of the following limsup set:
$$
S:= \limsup_{q\to \infty} \bigcup_{i\in \III(q)} B_i =
\bigcap_{Q=1}^\infty \bigcup_{q=Q}^\infty \bigcup_{i\in \III(q)}
B_i.
$$

{\bf Example.} The canonical example in this paper is the set of
simultaneously $\lambda$-well approximable points on a smooth
manifold $\MMM\subset \RR^n$, where $\lambda\ge 1/n$ is a fixed real
number. In this setting, $\III(q)$ denotes the set of points
$\vp/q\in \QQ^n$ for which there exists $\vx\in\MMM$ satisfying
\begin{equation}\label{eq7}
||q\vx-\vp||_\infty \le q^{-\lambda}.
\end{equation}
Then for each $i\in \III(q)$ we define the set $B^*_i$ by
$$
B^*_i:=\{\vx\in\MMM: \eqref{eq7} \mbox{ is satisfied }\}.
$$
Since $B^*_i$ is not necessarily a ball, we define $B_i$ to be a
ball of minimal diameter that contains $B^*_i$. If multiple such
balls exist, we choose any one of them. Note that for each
$i\in\III(q)$ we have the estimate $|B_i|\le
2\sqrt{n}q^{-1-\lambda}$. For some indices $i$ (specifically, when
$\vp/q$ lies close to $\MMM$) this upper bound may be nearly sharp
(i.e. $|B_i|\gg 2\sqrt{n}q^{-1-\lambda}$). For others, the diameter
of $B_i$ may be significantly smaller. In any case, we may take
$\rho(q) = 2\sqrt{d}q^{-1-\lambda}$.

We start with the classical result which is sometimes called the
Hausdorff-Cantelly lemma. For the proof, see for
example~\cite{ber_dod_1999}.

\begin{lemma}[Hausdorff-Cantelly]
If for given $s>0$ the following series converges
$$
\sum_{q=1}^\infty \sum_{i\in \III(q)} |B_i|^s
$$
then $\dim S\le s$.
\end{lemma}
Adapted to our setting, we can rewrite
\begin{corollary}\label{corl1}
$$
\dim S \le \inf\left\{s>0: \sum_{q=1}^\infty \#\III(q) \rho(q)^s
<\infty\right\}.
$$
\end{corollary}

This result provides a very good --- and often sharp --- upper bound
for $\dim S$, provided the quantity $\#\III(q)$ is well understood.
In the canonical example, if we take $\MMM=[0,1]^n$ then $\III(q)$
is the set of all rational points $\vp/q$ with $\vp\in [0,q]^n$ and
$\#\III(q) = (q+1)^n$. In this case, it follow easily from
Corollary~\ref{corl1} that $\dim S\le \frac{n+1}{\lambda+1}$ which
is a sharp bound due to the Jarnik-Besicovitch theorem. However, for
manifolds $\MMM$ of dimension strictly less than $n$, obtaining
precise estimates for $\# \III(q)$ is much more difficult. Moreover,
these values may fluctuate as $q$ varies.

\smallskip
{\bf Averaging over $q$.} We define the diadic blocks $\DDD(k):=
\{i\in \III(q): 2^k\le q<2^{k+1}\}$ and then rewrite
$$
S = \bigcap_{K=1}^\infty \bigcup_{k=K}^\infty \bigcup_{i\in \DDD(k)}
B_i.
$$
The idea is that, by grouping many sets $\III(q)$ into diadic blocks
$\DDD(k)$, we aim to smooth out local fluctuations in the sizes of
$\#\III(q)$. The quantity $\#\DDD(k)$ is often considerably easier
to estimate than $\#\III(q)$ itself. By applying the
Hausdorff-Cantelli lemma in this new setting we derive
\begin{corollary}\label{corl2}
Let $\rho^+(k):= \max\{\rho(q): 2^k\le q<2^{k+1}\}$. Then
$$
\dim S\le \inf\left\{s>0: \sum_{k=1}^\infty \#\DDD(k)
\rho^+(k)^s<\infty \right\}.
$$
Moreover, if $\rho(q)$ satisfies $\rho(q_1)\asymp \rho(q_2)$ as soon
as $0<q_1\le q_2<2q_1$ then this upper bound for $\dim S$ coincides
with that in Corollary~\ref{corl1}.
\end{corollary}

The proof of Corollary~\ref{corl2} is rather straightforward, so we
leave it as an exercise.

\medskip

{\bf Discrete partition.} One can split the sets $\DDD(k)$,
$k\in\ZZ_{\ge 0}$ into a finite number of subsets
\begin{equation}\label{eq8}
\DDD(k) = \bigcup_{j=1}^d \DDD_j(k).
\end{equation}
\begin{lemma}\label{lem2}
Let $d\in\NN$ be an absolute constant that does not depend on $q$.
Suppose that for each integer $k\ge 0$ the set $\DDD(k)$ is the
union~\eqref{eq8} of $d$ subsets. Then
$$
S = \bigcup_{j=1}^d S_j,
$$
where
$$
S_j:= \limsup_{k\to\infty} \bigcup_{i\in \DDD_j(k)} B_i.
$$
Therefore $\dim S\le \max\limits_{1\le j\le d} \dim S_j$.
\end{lemma}

\proof If $x\in S$ then it belongs to infinitely many sets
$\cup_{i\in \DDD(k)} B_i$. Then by Dirichlet principle, there exists
$1\le j\le d$ such that $x$ belongs to infinitely many sets
$\cup_{i\in \DDD_j(k)} B_i$ which is equivalent to saying that $x\in
S_j$. \endproof

The idea behind this approach is to partition each $\DDD(k)$ into
subsets $\DDD_j(k)$ in such a way that the cardinalities of these
smaller sets are easier to estimate. Lemma~\ref{lem2} is typically
applied when substantially different approaches are used to compute
upper bounds for $\dim S_j$ for each $1\le j\le d$.

\medskip
{\bf Continuous partition.} Another method is to introduce one or
more parameters for each set $\DDD(k)$ and then partition these sets
accordingly. Namely, consider the map
$$
\boldsymbol{\alpha}: \III\to \RR^d
$$
such that the image $\valp(\III)\subset \FFF$ for some bounded set
$\FFF\subset \RR^d$.

\begin{lemma}\label{lem3}
Suppose that for all $k\in\ZZ_{\ge 0}$, $\vp=(p_1,\ldots, p_d)\in
\FFF$ and all $\epsilon\in\RR^+$ one has
$$
\#\left(\valp^{-1}\left(\prod_{j=1}^d [p_j, p_j+\epsilon]\right)
\cap \DDD(k)\right) \le N(k,\vp,\epsilon),
$$
where $N(k,\vp,\epsilon): \NN\times \RR^d\times \RR^+\to \RR^+$.
Suppose also that the function
$$
M(k,\epsilon):= \frac {\log\left(\sup_{\vp\in \FFF}
N(k,\vp,\epsilon)\right)}{\log \rho^+(k)}
$$
has a limit as $\epsilon\to 0$ which is uniform in $k$. Then
\begin{equation}\label{lem3_eq}
\dim S\le \inf\left\{ s>0: \sum_{k=1}^\infty N(k)
\rho^+(k)^s<\infty\right\},
\end{equation}
where
$$
N(k) = \lim_{\epsilon\to 0}\, \sup_{\vp\in F}\, N(k,\vp,\epsilon).
$$
\end{lemma}

\proof Fix $\epsilon>0$ and divide $\FFF$ into $d$-dimensional
hypercubic regions of size $\epsilon$, i.e. each region is of the
form $\prod_{j=1}^d [p_j, p_j+\epsilon] \cap \FFF$. Since $\FFF$ is
bounded, there are only finitely many nonempty regions of that form
and their number does not depend on $k$ (but depends on $\epsilon$).
Let their number be $R$. Then we apply the discrete partition
technique and split each $\DDD(k)$ into $R$ subsets of the form
$$
\valp^{-1}\left(\prod_{j=1}^d [p_j, p_j+\epsilon]\right) \cap
\DDD(k)= : \DDD(k,\vp,\epsilon).
$$

By construction, we have $\#\DDD(k,\vp,\epsilon) \le
N(k,\vp,\epsilon) \le \sup_{\vp\in \FFF} N(k,\vp,\epsilon)=:
N(k,\epsilon)$. Then Corollary~\ref{corl2} implies that for
$S(\vp,\epsilon):= \limsup\limits_{k\to\infty} \bigcup_{i\in
\DDD(k,\vp,\epsilon)} B_i$,
$$
\dim S(\vp,\epsilon) \le \inf\left\{ s>0: \sum_{k=1}^\infty
N(k,\epsilon) \rho^+(k)^s<\infty\right\}.
$$
The series inside the infimum can be rewritten as
\begin{equation}\label{eq9}
\sum_{k=1}^\infty \rho^+(k)^{M(k,\epsilon)+s}.
\end{equation}

Now fix a small value $\delta>0$ and take $\epsilon$ such that for
all $k\in \ZZ_{\ge 0}$, $|M(k,\epsilon) - M(k,0)|\le \delta$. We can
do that since $M(k,\epsilon)$ has a uniform limit as $\epsilon\to
0$. Also notice that $\epsilon$ still does not depend on $k$, but
instead only depends on $\delta$. Therefore if the series
in~\eqref{lem3_eq} converges for the parameter $s$ then the
series~\eqref{eq9} converges for the parameter $s+\delta$ and then
Corollary~\ref{corl2} implies
$$
\dim S\le \max_{\vp\in \FFF} \dim S(\vp,\epsilon) \le \inf\left\{
s>0: \sum_{k=1}^\infty N(k) \rho^+(k)^s<\infty\right\} + \delta.
$$
By making $\delta$ arbitrarily small, the statement of the lemma
follows.
\endproof

\section{Outline of the result from~\cite{badziahin_2025} for $n=3$}

The beginning of the proof of Theorem~\ref{th1} is the same as
of~\cite[Theorem 2]{badziahin_2025} for $n=3$ where the same result
is proved but for a smaller range of $\lambda$. We outline the
required steps of the proof from there. The reader is encouraged to
see~\cite{badziahin_2025} for details.

From now on, $\MMM = \VVV_n(I):= \{(x,x^2,\ldots,x^n): x\in I\}$
where $I$ is any closed interval that does not contain zero. Next,
$\III = \{\vq \subset \NN\times \ZZ^n: \exists\, x\in\VVV_n\; \mbox{
s.t. \eqref{eq7} is satisfied}\}$. The value $q(i)$ for $i\in \III$
is just $q_0\in \NN$. Since $q_0$ plays a special role in the vector
$\vq$, we will sometimes use the notation $\vq = (q_0, \vq^+)$. The
ball $B_\vq$ is then the arc $\{(x,x^2,\ldots,x^n)\in \VVV_n:
\mbox{\eqref{eq7} is satisfied}\}$. Note that the projection map
$\pi:\VVV_n\to \RR$ to the first coordinate is bi-Lipschitz
therefore $\dim S = \dim \pi(S)$ and thus we can work with the set
of the first coordinates of $S$ instead of $S$ itself. The set $S$
is then $S_n(I, \lambda)$.

{\bf Step 1.} We apply an averaging approach. Then $\DDD(k)$ is equal
to $Q_n(I, \lambda, k)$ in terms of of~\cite{badziahin_2025}, which
is
$$
\DDD(k):= Q_n(I,\lambda,k) = \{\vq\in \III: 2^k\le q(\vq)<
2^{k+1}\}.
$$
For consistency, we will be using the notation $\DDD(k)$ in this
paper, but will provide its analogues from~\cite{badziahin_2025} for
easier referencing. Also, to simplify the notation, we denote $Q:=
2^k$ and notice that for all $\vq\in \DDD(k)$, $$\rho^+(k) =
\max_{\vq\in \DDD(k)}\{|B_\vq|\}\le 2\sqrt{n} Q^{-1-\lambda} \asymp
Q^{-1-\lambda}.$$

We introduce the following notation: we say that $a\gtrsim b$
(respectively, $a\lesssim b, a\simeq b$) if $Q^a \gg Q^b$
(respectively, $Q^a\ll Q^b, Q^a\asymp Q^b$).

{\bf Step 2.} Split $I$ into several intervals of the
form $B(x_m,Q^{-\frac{1+\lambda}{2}})$ (or rather of radius $\asymp
Q^{-\frac{1+\lambda}{2}}$ if $|I|$ is not an integer multiple of
that number). Then $\DDD(k)$ splits into subsets $\DDD(k,m)$ of the form
\begin{equation}\label{eq10}
Q_n(I,\lambda,k,m):= \left\{\vq\in Q_n(I,\lambda,k): \begin{array}{l}
|q_0|\ll Q;\\
|q_0x_m - q_1|\ll Q^{\frac{1-\lambda}{2}};\\
|(1-i)x_m^iq_0 + ix_m^{i-1} q_1 - q_i|\ll Q^{-\lambda},\; 2\le
i\le n.
\end{array}\right\}
\end{equation}
The idea is that for $x\in B(x_m, Q^{-\frac{1+\lambda}{2}})$ the
piece of the curve $\VVV_n$ can be treated as a straight segment,
and therefore all rational points close to it must lie inside a
convex box that is defined by~\eqref{eq10}. We denote such a box by
$\Delta_m$.

{\bf Step 3.} Apply discrete partitioning to $\DDD(k)$. We write
$\DDD(k) = \DDD_1(k)\sqcup \DDD_2(k)$, where $\DDD_1(k) =
Q_n^1(I,\lambda,k)$ consists of the union of all  $\DDD(k,m)$ such
that $\#\DDD(k,m) \le Q^{\frac{3-(2n-1)\lambda}{2}}$. Then we
show~\cite[Lemma 2]{badziahin_2025} that $\dim S_1 \le
\frac{2-(n-1)\lambda}{1+\lambda}$, therefore it only remains to
estimate $\dim S_2$.

One can check that the volume of $\Delta_m$ is
$Q^{\frac{3-(2n-1)\lambda}{2}}$. If its last successive minimum
satisfies $\tau_{n+1}\le 1$ then $\Delta_m$ contains $n+1$ linearly
independent vectors which implies
$$
\# \DDD(k,m) \asymp \vol(\Delta_m) \ll Q^{\frac{3-(2n-1)\lambda}{2}}\quad \Longrightarrow\quad \DDD(k,m)\subset \DDD_1(k).
$$
The idea here is that for a ``generic'' convex centrally symmetric
figure all successive minima are of approximately the same size,
i.e. $\tau_1\asymp \tau_2\asymp\cdots\asymp \tau_{n+1}$. Therefore,
in view of the second Minkowski theorem,
$\tau_1\tau_2\cdots\tau_{n+1} \asymp Vol^{-1}(\Delta_m)$, one has
for $\lambda\le \frac{3}{2n-1}$ that a ``generic'' box $\Delta_m$
has $\tau_{n+1}<1$ and therefore the points from the ``majority'' of
boxes $\Delta_m$ belong to $\DDD_1(k)$. So now we are left with the
set $\DDD_2(k)$ from the ``exceptional'' boxes.

{\bf Step 4.} Apply continuous partitioning to $\DDD_2(k)$. For each
of the remaining sets $\DDD(k,m)\subset\DDD(k)$ we associate the
parameter $\delta=\delta_m$ which is given by $\tau_{n+1} =
Q^{\delta}$. By construction, we always have $\delta\ge 0$. On the
other hand, since $CQ^\lambda \Delta_m$ for appropriately chosen
absolute constant $C$ always contains the unit cube centered at $0$,
we have $\delta\le \lambda$. Notice that all the points $\vq\in
\DDD(k,m)$ lie in some proper linear subspace of $\RR^{n+1}$. We
denote the hyperplane of the smallest height that contains all
$\vq\in \DDD(k,m)$ by $\PPP(m)$ and denote its equation by
$\va(m)\cdot \vq = 0$. In other words, for $\vq\in \DDD(k,m)$ we
associate the vector $\va = \va(m)$. Observe that for all $\vx\in
B_\vq$ one has $\va\cdot \vx = a_0 + a_1x+\cdots + a_nx^n =
P_\va(x)$ and
\begin{equation}\label{eq11}
|P_\va(x)| = |q_0^{-1} (\va\cdot\vq + a_1(q_0x-q_1) +
\cdots + a_n(q_0x^n-q_n)| \ll ||\va||_\infty
Q^{-1-\lambda}.
\end{equation}

Next, for each $\vq\in \DDD(k,m) \cap \DDD_2(k)$  we also associate
the interval $J(m)$ such that it contains the balls $B_\vq$ for all
$\vq\in \DDD(k,m)$; has the maximal possible length and for all
$x\in J(m)$ the inequality~\eqref{eq11} is satisfied. Sometimes we
will denote this interval by $J(\vq)$ to emphasize that it is
associated with a particular vector $\vq\in \DDD_2(k)$. We introduce
the parameter $\eta = \eta_m\in \RR$ as follows:
\begin{equation}\label{eq14}
|J(m)| = Q^{-\frac{1+\lambda}{2} - \eta}.
\end{equation}
The parameter $\eta$ can be negative as well as positive. If
$\eta<0$ then the same interval $J$ can be associated with several
consecutive sets $\DDD(k,m)$ (in fact, up to $Q^\eta$ of
them). Since $J(m)$ contains at least one of the balls $B_\vq$ we
have $\eta\le \frac{1+\lambda}{2}$. On the other hand, we also
always have $|J(m)|\ll 1$ therefore $\eta\ge -\frac{1+\lambda}{2}$.
The upshot is that the set of parameters $\vp:= (\delta, \eta)$ lies in
a bounded region $\FFF$.

Now, we apply continuous partitioning to the sets $\DDD_2(k)$
with respect to the pair of parameters $\vp=(\delta, \eta)$. Then for a
given $\vp\in\FFF$, all $\vq\in \DDD_2(k,\vp,\epsilon)$
satisfy (see~(20) in~\cite{badziahin_2025})
\begin{equation}\label{eq12}
||\va||_\infty \ll \left\{\begin{array}{ll}
Q^{\lambda-\eta-\delta+\epsilon}& \mbox{if }\; \eta\ge 0;\\[1ex]
Q^{\lambda-\delta+\epsilon}& \mbox{if }\; \eta<0.
\end{array}\right.
\end{equation}
Also, the number of vectors $\vq\in \DDD_2(k,\vp,\epsilon)$ that share the
same interval $J(m)$ is bounded from above by (see~(21) in
\cite{badziahin_2025})
\begin{equation}\label{eq13}
\ll\left\{ \begin{array}{ll} Q^{\frac{3-5\lambda}{2} +
2\delta+2\epsilon}&\mbox{if }\; \eta\ge 0;\\[1ex]
Q^{\frac{3-5\lambda}{2} + 2\delta -\eta+3\epsilon}&\mbox{if }\;
\eta<0.
\end{array}\right.
\end{equation}
In order to estimate $\#\DDD_2(k,\vp,\epsilon)$ it remains to bound
from above the number of polynomials $P_\va(x)$ of degree $n$ whose
height $\va$ satisfies~\eqref{eq12} and whose values
satisfy~\eqref{eq11} for all $x$ in the interval $J$ of length at
least as in~\eqref{eq14}. Denote the set of such polynomials by
$\AAA_n(I,\lambda, k, \vp,\epsilon)$. If there is no confusion about
its set of parameters, we will call it $\AAA(k)$ or $\AAA(k,
\vp,\epsilon)$.

{\bf Step 5.} Discretely partition the set $\DDD_2(k) = \DDD_3(k)
\cup \DDD_4(k)$, where for $\vq\in \DDD_3(k)$ the discriminant of
the corresponding polynomial $P_\va$ is non-zero.  The polynomials
$P_\va$ for $\vq\in \DDD_4(k)$ then have zero discriminant.
Consequently, the corresponding set $\AAA(k)$ also splits in two:
$\AAA_3(k)\cup \AAA_4(k)$.

{\bf Step 6.} Narrow the collection $\AAA_{3,n}(I,\lambda,k,
\vp,\epsilon)$ of polynomials $P_\va(x)$. Let $C = (n+1)^{n+3}2^n
\asymp 1$ and $j\in\{0,\ldots, n\}$. Define $R_j(x):= (Cx)^3
P_\va((Cx)^{-1} + j)$ and $f_j(x) = \frac{1}{C(x-j)}$. we
show~\cite[Lemma 6]{badziahin_2025} that for at least one of $j$ in
this range, the polynomial $R_j(x)$ satisfies~\eqref{eq12}
and~\eqref{eq11} for all $x\in f_j(J)$ whose length is $\gg
Q^{-\frac{1+\lambda}{2}-\eta}$ which is comparable to~\eqref{eq14}.
On top of that, the leading coefficient of $R_j$ has the largest
absolute value among all the coefficients of $R_j$. Finally, it is
not hard to verify that the discriminant of $R_j$ is not zero. This
implies
$$
\#\AAA_{3,n}(I,\lambda,k, \vp,\epsilon) \ll \sum_{j=0}^3
\#\AAA_{3,n}^{*}(I_j, \lambda, k,\vp,\epsilon),
$$
where $\AAA_{3}^{*}(k,\vp,\epsilon)$ consists of those $\va\in
\AAA_3(k,\vp,\epsilon)$ that satisfy $||\va||_\infty = |a_n|$.
$I_j\subset \frac{1}{C(I-j)}$ are some intervals that are separated
from zero.

Therefore, it is sufficient to estimate the cardinality of
$\AAA_n^{*} (I, \lambda,k)$ for all intervals $I$ distanced from
zero. In the further discussion we will focus on these sets and for
convenience lift the star from the superscript.

\section{Sketch of the new development}

In the further discussion we will focus on the case of $n=3$.
However, wherever possible, we will make proofs for arbitrary values
of $n$, as this may help with further research. The proof of
Theorem~\ref{th1} involves many steps and cases. To help navigate
through all of them we first provide a schematic outline of the
proof and then will go through each of its steps.

\begin{itemize}
%\begin{itemize}
\item[\bf 7.1] Apply continuous partitioning to $\DDD_3(k)$ by introducing the parameter $\kappa$ such that the distance between the furthest zeroes of $P_\va$ is $Q^{-\kappa}$.
\item[\bf 7.2] Prove Theorem~\ref{th3} and apply it to estimate $\#\DDD_3(k,\vp, \epsilon)$ which in turn will imply Theorem~\ref{th1} for the set $S_3$.
%\end{itemize}

%\begin{itemize}
\item[\bf 8.1] Rule out the case of $\delta\ge 2\lambda-1$ for $\DDD_4(k)$.
\item[\bf 8.2] Discretely partition $\DDD_4(k) = \DDD_5(k)\cup \DDD_6(k)\cup \DDD_7(k)$,
 where rational numbers $q_1/q_0$ for $\vq\in \DDD_5(k)$ and $\vq\in \DDD_6(k)$
 lie close enough to another rational number $u/v$ with a very small denominator.
 Namely, conditions~\eqref{eq26} and~\eqref{eq27} are satisfied for $\DDD_5(k)$
 and $\DDD_6(k)$ respectively.
\item[\bf 8.3] Establish Theorem~\ref{th1} for $S_5$ and $S_6$.
\item[\bf 8.4] Write $P_\va(x) = (ax-b)^2(cx-d)$. Discretely partition $\DDD_7(k) = \DDD_8(k)\cup\DDD_9(k)\cup\DDD_{10}(k)$ in the following way. For $\vq\in \DDD_8(k)$, the interval $J(\vq)$ contains $d/c$ but not $b/a$. For $\vq\in \DDD_9(k)$, $J(\vq)$ contains $b/a$ but not $d/c$. Finally, for $\vq\in \DDD_{10}(k)$, there exists a rational number $u/v\in J(\vq)$ with $v\le H^{1/3}$.
\item[\bf 8.5] Establish Theorem~\ref{th1} for $S_8$.
\item[\bf 8.6] Continuously partition the set $\DDD_9(k)$ by introducing the parameter $\sigma$ such that $|c| = Q^\sigma$. Then verify Theorem~\ref{th1} for $S_9$.
\item[\bf 8.7.1] Continuously partition the set $\DDD_{10}(k)$ by introducing the parameter $\delta^*$ such that the third successive minimum of the corresponding $\Delta_m$ is $\tau_3= Q^{\delta^*}$.
\item[\bf 8.7.2] Rule out the case $\delta^* + \delta \le \frac12(1-\lambda) + \epsilon$ for $\DDD_{10}(k)$.
\item[\bf 8.7.3] Discretely partition $\DDD_{10}(k) = \DDD_{11}(k)\cup\DDD_{12}(k)$ where for all  points $\vq\in\DDD_{11}(k)$ that lie in the same $\Delta_m$ the condition~\eqref{eq28} is satisfied. The set $\DDD_{12}(k)$ then contains all the remaining points. Verify Theorem~\ref{th1} for $S_{11}$.
\item[\bf 8.7.4] Show that $\DDD_{12}(k)$ is empty and hence establish Theorem~\ref{th1} for this set.
%\end{itemize}
\end{itemize}

\section{The case of non-zero discriminant}

Given a polynomial $P\in \CC[x]$ of degree at least 2, we define the
notion $r(P)$ to be the largest distance between its roots. We will
focus on the set $\AAA_3(k)$ and modify Lemma~7
from~\cite{badziahin_2025} by introducing the parameter $\kappa$
given by $r(P_\va)=Q^{-\kappa}$. In other words, we will now deal
with the triples of parameters $\vp=(\delta, \eta, \kappa)$. The
idea is that if the determinant $D(P_\va)$ of the polynomial is much
less than $a_n^{2n-2}$ then its parameter $r(P_\va)$ is most likely
much less than 1. This fact will later help us better estimate the
number of suitable polynomials $P_\va$.

By Cauchy's bound on the polynomial roots, as soon as
$H(P_\va)=a_n$, we always have $r(P_\va)\ll 1$, i.e. $\kappa\gtrsim
0$. On the other hand,
$$
1\le |D(P_\va)| \le a_n^{2n-2}Q^{-n(n-1)\kappa}
$$
which implies $Q^\kappa \le a_n^{2/n}$. In view of~\eqref{eq12}, we
derive that $\kappa$ varies within a bounded set.

\begin{lemma}\label{lem6}
Suppose that for a given polynomial $P_\va$ with $||\va||_\infty =
a_n$ there exist $w,\kappa,\eta\in \RR$ with $w>0$, $-\frac{w}{2} <
\eta < \frac{w}{2}$, and the interval $J$ of length $|J| \gg
Q^{-\frac{w}{2}-\eta}$ such that $r(P_\va) = Q^{-\kappa}$ and
$\forall x\in J$, $|P_\va(x)| < a_nQ^{-w}$. Then the discriminant of
$P_\va$ satisfies
\begin{equation}\label{lem1_eq}
D(P_\va) \ll a_n^{2n-2} Q^{-w-(n-1)(n-2)\kappa+2\eta}.
\end{equation}
\end{lemma}

{\bf Remark.} Substituting the bound $\kappa\gtrsim 0$
into~\eqref{lem1_eq}, gives the initial bound
from~\cite[Lemma~7]{badziahin_2025}.

While the proof is similar to that in~\cite{badziahin_2025}, we
reproduce it here for convenience of readers.

\proof Fix a point $x_0\in J$ and consider any $x\in J$. We get
$$
P_\va(x) = P_\va(x_0) + (x-x_0)P'_\va(x_0) + \ldots +
\frac{1}{n!}(x-x_0)^n P_\va^{(n)}(x_0).
$$
Let $y_1, y_2, \ldots, y_{n+1}\in J$ be such that $y_1$ and
$y_{n+1}$ are the endpoints of $J$ and $y_2-y_1 = \ldots =
y_{n+1}-y_n$. That immediately implies $|y_{i+1}-y_i|= |J|/n \asymp
|J|$ for all $1\le i\le n$. Also denote $b_i =
\frac{P_\va^{(i)}(x_0)}{i!}$. Then the values $b_i$ are the
solutions of the following matrix equation
$$
\left(\begin{array}{ccccc} 1&y_1-x_0&(y_1-x_0)^2&\cdots&
(y_1-x_0)^n\\
1&y_2-x_0&(y_2-x_0)^2&\cdots&
(y_2-x_0)^n\\
\vdots&\vdots&\vdots&\ddots&\vdots\\
1&y_{n+1}-x_0&\cdots&\cdots&(y_{n+1}-x_0)^n
\end{array}\right) \left(\begin{array}{c}
b_0\\b_1\\\vdots\\b_{n}
\end{array}\right) = \left(\begin{array}{c}
P_\va(y_1)\\P_\va(y_2)\\\vdots\\P_\va(y_{n+1})
\end{array}\right).
$$
Notice that on the left hand side we have the Vandermonde matrix.
Let's call it $V$. Since $|y_j-y_i|\gg Q^{-\frac{w}{2} - \eta}$ for
all $1\le i<j\le n+1$, its determinant is
$$
\det V \gg \left(Q^{-\frac{w}{2}-\eta}\right)^{\frac{n(n+1)}{2}}.
$$
Then Cramer's rule gives for $2\le i\le n$
\begin{equation}\label{eq1}
|P_\va(x_0)| \asymp |b_0| \ll a_nQ^{-w};\; |P'_\va(x_0)| \asymp
|b_1| \ll a_nQ^{-\frac{w}{2}+\eta};\; |P^{(i)}_\va(x_0)|\ll
a_nQ^{-\frac{2-i}{2}w + i\eta}.
\end{equation}

Let $x_1,x_2,\ldots, x_n$ be the roots of $P_\va$ such that
$|x_1-x_0|\le |x_2-x_0|\le\cdots\le |x_n-x_0|$. Notice that for all
$x$ in the segment between $x_1$ and $x_0$  and all $2\le i\le n$ we
get $|x-x_i|\ll |x_0-x_i|$. Together with~\eqref{eq1}, that implies

$$
a_nQ^{-\frac{w}{2}+\eta}\gg|P'_\va(x_0)|\gg |P'_\va (x_1)| = a_n
|(x_1-x_2)\cdots (x_1-x_n)|.
$$
For all other distances $|x_i-x_j|, 2\le i<j\le n$, we use the bound
$|x_i-x_j|\le Q^{-\kappa}$. Combining all these bounds together
gives
$$
|D(P_\va)| = a_n^{2n-2} \prod_{1\le i<j\le n} |x_i-x_j|^2 \ll
a_n^{2n-2} Q^{-w+2\eta -(n-1)(n-2)\kappa}.
$$
\endproof

Let $P_\va\in \ZZ[x]$ be a given cubic polynomial and $R,H>0$
positive real numbers. Denote by $N(P_\va,H,R)$ the number of
polynomials $P$ in the same equivalence class as $P_\va$ such that
$H(P)\le H$ and $r(P)\ge R^{-1}$.

\begin{proposition}\label{prop1}
For any $\epsilon>0$ there exists $c=c(\epsilon)>0$ such that for
any cubic polynomial $P_\va$ one has
$$
N(P_\va,H,R)\le c (\min\{H^{2/3+\epsilon}D(P_\va)^{-1/6}, \log
H\cdot R\}+ H^\epsilon).
$$
\end{proposition}

\proof For a given polynomial $P(x) = c_0 + c_1x + c_2x^2 + c_3x^3$
we introduce the following height
$$
H_d(P):= \max\{|c_2|, |c_3|, |c_1c_2|^{1/2}, |c_0c_2^3|^{1/4},
|c_0c_3|^{1/2}, |c_1^3c_3|^{1/4}, |c_0c_1c_2c_3|^{1/4}\}.
$$
Let $R_\va$ be the polynomial $c_0+c_1x+c_2x^2+c_3x^3$ with the
minimal height $H_d$ among all polynomials in the equivalence class
of $P_\va$. Let $x_1,x_2,x_3$ be its roots. In~\cite[Lemma
9]{badziahin_2025} it is shown that $|x_i-x_j|\gg 1$ for all $1\le
i<j\le 3$. Order the roots in such a way that $|x_1-x_2|$ is the
smallest distance and $|x_2-x_3|$ is the largest one among the roots
of $R_\va$. Denote $d\ge 0, D\ge 0$ in such a way that $2^d\le
|x_1-x_2|<2^{d+1}$ if $|x_1-x_2|\ge 1$ and $d=0$ otherwise. The
parameter $D$ is analogously defined for the distance $|x_2-x_3|$.
Then we have
$$
|x_1-x_2|\asymp 2^d, |x_1-x_3|\asymp |x_2-x_3|\asymp 2^D.
$$
This in turn implies $|D(P_\va)| = |D(R_\va)| \asymp c_3^4
2^{2d+4D}$.

We need to compute an upper bound for the number of M\"obius
transforms $\mu$ such that $H((c-ax)^3R_\va\circ \mu^{-1})\le H$ and
$r((c-ax)^3R_\va\circ \mu^{-1})\ge R^{-1}$. Up to an absolute
constant, this number is bounded by the number of pairs $a,b$ such
that $\mu(x) = \frac{cx+d}{ax+b}$, $|a|\ge |c|, |b|\ge |d|$ and the
resulting polynomial $P = (c-ax)^3R_\va\circ \mu^{-1}$ has bounded
height and $r(P)$. Denote the set of such pairs by $M(P_\va,H,R) =
M$. One can check that the leading coefficient $c_3(\mu)$ of $P$
equals (see~\cite[Equation (31)]{badziahin_2025})
\begin{equation}\label{eq2}
c_3(\mu) = a^3R_\va(-b/a) = -c_0a^3+c_1a^2b-c_2ab^2+c_3b^3.
\end{equation}

Fix $1>\epsilon>0$, $2> t\ge 0$ and consider the set
$\SSS(t,\epsilon)$ of points $a,b\in \ZZ^2$ such that
$$
|a|^{-t-\epsilon}<\left|y_1 + \frac{b}{a}\right|\le |a|^{-t}
$$
where $y_1$ is one of the roots of $R_\va$ closest to $-b/a$. Let
$y_2,y_3$ be the remaining roots of $R_\va$. One can check that the
number of pairs in this set such that $|a|\le A$ is bounded by $\ll
A^{2-t}$.

Consider $(a,b)\in \SSS(t,\epsilon)\cap M$. Since all the roots of
$R_\va$ are placed far apart from each other, we must have $|y_j +
b/a|\asymp|y_j-y_1|\gg 1$, $j\in\{2,3\}$. Then from~\eqref{eq2} we
compute
$$
|c_3(\mu)| = |a^3R_\va(-b/a)| \asymp |a^3c_3
(y_1-y_2)(y_1-y_3)(y_2-y_3)|\cdot \frac{|y_1+b/a|}{|y_2-y_3|} \gg
\left| \frac{a^{3-t-\epsilon}|D(P_\va)|^{1/2}}{c_3(y_2-y_3)}\right|.
$$
Since we must have $|c_3(\mu)|\le H$, this establishes an upper
bound on the size $|a|$:
\begin{equation}\label{eq4}
|a|\ll \left(\frac{H
|c_3(y_2-y_3)|}{|D(P_\va)|^{1/2}}\right)^{\frac{1}{3-t-\epsilon}}\ll
\left(\frac{H}{|2^{d/2}D(P_\va)|^{1/4}}\right)^{\frac{1}{3-t-\epsilon}}
.
\end{equation}

Now we estimate the value of $r(P)$. By construction, the roots of
$P$ are $\mu(y_1), \mu(y_2)$ and $\mu(y_3)$. We compute
$$
|\mu(y_1)- \mu(y_i)| = \frac{|y_1-y_i|}{|(ay_1+b)(ay_i+b)|} \gg
\frac{|y_1-y_i|}{|a|^{2-t} |y_1-y_i|}= |a|^{-2+t},
$$
where $i\in\{2,3\}$ and
$$
|\mu(y_2)-\mu(y_3)| = \frac{|y_2-y_3|}{|(ay_2+b)(ay_3+b)|} \asymp
\frac{|y_2-y_3|}{a^2 |y_1-y_3|\cdot |y_1-y_2|}\ll a^{-2}.
$$
From the last two bounds we see that $r(P)\asymp \max_{2\le i\le
3}|\mu(y_1) - \mu(y_i)|$ and then
$$
r(P) \ll \max_{2\le i\le
3}\left\{\frac{|y_1-y_i|}{|(ay_1+b)(ay_i+b)|}\right\}\ll |a|^{-2+t}.
$$
This leads to the inequality
\begin{equation}\label{eq3}
|a|\ll R^{\frac{1}{2-t}}.
\end{equation}

Combining~\eqref{eq2} and~\eqref{eq3}, gives
\begin{equation}\label{eq29}
\# \SSS(t,\epsilon)\cap M \ll
\min\left\{(HD(P_\va)^{-1/4})^{\frac{2-t}{3-t} + \epsilon_1},
R\right\}\le \min\{H^{\frac23+\epsilon_1}|D(P_\va)|^{-\frac16},R\}
\end{equation}
where $\epsilon_1>0$ tends to zero together with $\epsilon$.

Next, consider the set $\SSS(2)$ of pairs $(a,b)$ such that
$|y_1+\frac{b}{a}|\le a^{-2}$. Then one can check that
$$
\# \SSS(2)\cap M \ll (HD(P_\va)^{-1/4})^{\epsilon_1}.
$$
This is shown in~\cite[Proposition 5, cases $S(2,1/2)$ and
$S(5/2)$]{badziahin_2025}.

The remaining pairs $a,b$ satisfy $|y_i + a/b|>1$ for all
$i\in\{1,2,3\}$. Split them into subsets $\SSS_0(k)$ where every
$(a,b)\in \SSS_0(k)$ satisfies
$$
2^k\le \left|y_1 + \frac{b}{a}\right|<2^{k+1}.
$$
The number of pairs with $|a|<A$ in every such set equals $\asymp2^k
A^2$. First consider $k\le d$. In this case the analogous inequality
to~\eqref{eq4} for $a\in \SSS_0(k)\cap M$ is
$$
|a|\ll
\left(\frac{H}{2^{k+d/2}|D(P_\va)|^{1/4}}\right)^{\frac{1}{3}} .
$$
This implies
$$
\# (\SSS_0(k)\cap M) \ll 2^{\frac{k-d}{3}}H^{\frac23}
|D(P_\va)|^{-1/6}.
$$

Analogously to the case of $\SSS(t,\epsilon)$, we compute that
$|\mu(y_2)-\mu(y_3)|\ll a^{-2}2^{-d}$ and $$|\mu(y_1)-\mu(y_i)|
\asymp (2^ka^2)^{-1} \gg |\mu(y_2)-\mu(y_3)|.$$ Therefore $r(P_\va)
\ll (2^k a^2)^{-1}$ which implies $2^k a^2\ll R$ and hence
\begin{equation}\label{eq5}
\#\SSS_0(k)\cap M \ll R.
\end{equation}

Next, let $d<k\le D$. If $|y_1-y_2|\asymp 2^d$ or $|y_1-y_3|\asymp
2^d$ then
$$
|c_3(\mu)|=\left|a^3
c_3\left(y_1+\frac{b}{a}\right)\left(y_2+\frac{b}{a}\right)\left(y_3+\frac{b}{a}\right)\right|\asymp
|a^3c_32^{2k+D}| \asymp |a^3 D(P_\va)^{1/4} 2^{2k-d/2}|
$$
and
\begin{equation}\label{eq6}
|a|\ll
\left(\frac{H}{2^{2k-d/2}|D(P_\va)|^{1/4}}\right)^{1/3}\qquad\Longrightarrow\qquad
\#(\SSS_0(k)\cap M)\ll 2^{\frac{d-k}{3}}
H^{\frac23}|D(P_\va)|^{-1/6}.
\end{equation}

Without loss of generality assume that $|y_1-y_2|\asymp 2^d$. then
one computes $|\mu(y_1)-\mu(y_2)|\asymp (2^{2k-d}a^2)^{-1}$ and
$|\mu(y_i)-\mu(y_3)|\asymp (2^ka^2)^{-1}$ for $i\in\{1,2\}$.
Therefore $r(P_\va)\asymp (2^ka^2)^{-1}$ and~\eqref{eq5} follows
again.

If $|y_2-y_3|\asymp 2^d$ then analogous computations give
$$
|c_3(\mu)|\asymp |a^3c_32^{2D+k}| \asymp |a^3 D(P_\va)^{1/4} 2^{k +
D -d/2}|\ge  |a^3 D(P_\va)^{1/4} 2^{2k-d/2}|
$$
and the same bound~\eqref{eq6} for $\#(\SSS_0(k)\cap M)$ holds.
Also, we compute $|\mu(y_1)-\mu(y_i)|\asymp (2^ka^2)^{-1}$,
$|\mu(y_2)-\mu(y_3)|\asymp (2^{2D-d}a^2)^{-1}\ll
|\mu(y_1)-\mu(y_i)|$ and hence~\eqref{eq5} is satisfied.

Finally, for $k>D$, $|c_3(\mu)|\asymp |a^3 D(P_\va)^{1/4} 2^{3k -
D-d/2}|\ge |a^3 D(P_\va)^{1/4} 2^{2k-d/2}|$ and hence we
get~\eqref{eq6}. In this case we also have $|\mu(x_i)-\mu(x_j)|\ll
(2^ka^2)^{-1}$ for all $1\le i<j\le 3$ and hence~\eqref{eq5} is
satisfied.

Notice that from~\eqref{eq6} we get that for $k>2\log_2 H + C$ for
large enough constant $C$ the set $\SSS_0(k)\cap M$ is empty.
Therefore the number of $k$ for which it is nonempty, is $\ll \log
H$.

To finish the proof of the proposition, we split the interval
$[0,2)$ into $N$ subintervals of equal length $\epsilon$. Then we
split the set $M(P_\va,H,R)$ into subsets
$$
\bigcup_{i=0}^N (\SSS(i\epsilon, \epsilon)\cap M) \bigcup
(\SSS(2)\cap M)\bigcup_{k=0}^\infty (\SSS_0(k)\cap M) .
$$
By combining the estimates~\eqref{eq29},~\eqref{eq5}
and~\eqref{eq6}, the total number of points in this union is bounded
from above by
$$
\min\left\{\left(N + 2\sum_{k=0}^\infty 2^{-k/3}\right)H^{\frac23 +
\epsilon_1} D(P_\va)^{-1/6}, (N+1+\log H)R\right\} +
(HD(P_\va)^{-1/4})^{\epsilon_1}.
$$
Here $\epsilon_1$ can be taken arbitrary small and $N =
2\epsilon^{-1}$. The conclusion of the proposition then follows
immediately. \endproof

\begin{theorem}\label{th3}
For any $\epsilon>0$ there exists a constant $c = c(\epsilon)$ such
that the number $N(H,D,R)$ of polynomials $P$ with $H(P)\le H$,
$0<|D(P)|\le D$ and $r(P)\ge R^{-1}$ is bounded from above by
\begin{equation}\label{th3_eq}
N(H,D,R)\le c (\min\{H^{2/3+\epsilon} D^{5/6}, \log H\cdot DR\} +
H^\epsilon D).
\end{equation}
\end{theorem}

\proof Let $h(d)$ be the number of equivalence classes of cubic
polynomials that share the discriminant $d\neq 0$. For convenience
of notation we set $h(0) = 0$. Davenport~\cite{davenport_1951}
showed that
$$
\sum_{d=-D}^D h(d) \asymp D.
$$

By Proposition~\ref{prop1}, for any equivalence class of cubic
polynomials of discriminant $d$ we have at most
$$
\ll \min\{H^{2/3+\epsilon}d^{-1/6}, \log H\cdot R\} + H^\epsilon
$$
polynomials $P$ with $H(P)\le H$ and $r(P)\ge R^{-1}$. Summing over
all such classes gives
$$
N(H,D,R) \ll \sum_{d=-D}^D (\min\{H^{2/3+\epsilon}d^{-1/6}, \log
H\cdot R^{1+\epsilon}\} + H^\epsilon)
$$
The application of Abel's summation formula finishes the proof
\endproof

Now we are ready to estimate $\#\AAA_3(k,\vp,\epsilon)$ and thus
conclude the step~7.2 of the proof of Theorem~\ref{th1}. Recall that
for each $\va \in \AAA_3(k,\vp,\epsilon)$ the corresponding
polynomial $P_\va$ satisfies the following conditions: $H(P_\va)$ is
bounded by~\eqref{eq12}, $r(P_\va)\ge Q^{-\kappa-\epsilon}$ and in
view of Lemma~\ref{lem6},
$$
D(P_\va)\ll H^4 Q^{-1-\lambda - 2\kappa + 2\eta+2\epsilon},
$$
where for simplicity we denote by $H$ the upper bound on $H(P_\va)$.
Theorem~\ref{th3} can now be applied. To simplify the computations,
we bound the second term of the minimum by $H^\epsilon DR$. In this
case we get
$$
H^\epsilon D\le H^\epsilon DR\quad\mbox{and}\quad H^\epsilon D\le
H^{2/3+\epsilon}D^{5/6}.
$$
Therefore the term $H^\epsilon D$ can be ignored. Also we only use
the second term in the minimum from~\eqref{th3_eq}. Then we derive
$$
\# \AAA_3(k,\vp,\epsilon)\ll H^{4+\epsilon}Q^{-1-\lambda -2\kappa +
2\eta + \kappa+ 3\epsilon}.
$$
For $\eta\ge 0$ this bound together with~\eqref{eq12}
and~\eqref{eq13} imply
$$
\#\DDD_3(k,\vp,\epsilon) \ll Q^{\frac{3-5\lambda}{2} +
2\delta+2\epsilon +
(4+\epsilon)(\lambda-\eta-\delta+\epsilon)-1-\lambda-\kappa+2\eta+3\epsilon}
\ll Q^{\frac{1+\lambda}{2} - \kappa-2\eta-2\delta + c_1\epsilon}$$
for some absolute constant $c_1>0$. Clearly, the right hand side
attains its maximum for $\eta=\delta=\kappa=0$.

For $\eta<0$ analogous computations reveal
$$
\#\DDD_3(k,\vp,\epsilon) \ll Q^{\frac{3-5\lambda}{2} +
2\delta-\eta+3\epsilon +
(4+\epsilon)(\lambda-\delta+\epsilon)-1-\lambda-\kappa+2\eta+3\epsilon}
\ll Q^{\frac{1+\lambda}{2} - \kappa+\eta-2\delta + c_2\epsilon}.
$$
Here again, the right hand side maximises when
$\delta=\eta=\kappa=0$. As an upshot, one can choose the value of
$N(k,\vp,\epsilon)$ from Lemma~\ref{lem3} to be such that
$\sup\limits_{\vp\in\FFF} N(k,\vp,\epsilon) = Q^{\frac{1+\lambda}{2}
+ \max\{c_1,c_2\}\epsilon}$.

Compute
$$
M(k,\epsilon) = \frac{\log (\sup_{\vp\in \FFF}
\#\DDD_3(k,\vp,\epsilon))}{\log \rho^+(k)} = \frac{(1+\lambda)/2 +
\max\{c_1,c_2\}\epsilon}{1+\lambda}
$$
This expression uniformly tends to $\frac12$ as $\epsilon\to 0$.
Therefore Lemma~\ref{lem3} states that $\dim S_3 \le \frac12$. One
can easily check that for $\lambda\le \frac35$ this is smaller than
$\frac{2-2\lambda}{1+\lambda}$.

\section{Auxiliary results for the case of zero discriminant: the proof of steps 8.1 -- 8.3}

From now on, we consider the set $\DDD_4(k)$ and the corresponding
set of hyperplanes (or equivalently, polynomials) $\AAA_4(k)$.
Recall that for each $\va\in \AAA_4(k)$ the polynomial $P_\va(x)$
has zero discriminant, i.e. it is of the form $P_\va(x) =
(ax-b)^2(cx-d) = P_1^2P_2$ for some integer $a,b,c,d$.

Fix $\va = \va_m\in A_4(k)$ and consider the corresponding interval
$J = J(m)$. Then for all $x\in J$ we have
\begin{equation}\label{eq15}
|P_\va(x)| \stackrel{\eqref{eq11}}\ll H(P_\va)Q^{-1-\lambda}
\stackrel{\eqref{eq12}}\ll
H(P_\va)^{1-\frac{1+\lambda}{\lambda-\delta+ \epsilon}}.
\end{equation}
In view of the Gelfond Lemma we have $H(P_\va) \asymp
H^2(P_1)H(P_2)$, therefore there exists $i\in\{1,2\}$ such that
$$
|P_i(x)|\ll H(P_i)^{1-\frac{1+\lambda}{\lambda-\delta+\epsilon}}.
$$
If a transcendental $x$ belongs to $S_4(\vp,\epsilon)$, i.e. it
satisfies infinitely many inequalities~\eqref{eq15}, then it must
satisfy infinitely many above inequalities for linear polynomials
$P$. By the Jarnik-Besicovich theorem we then have
$$
\dim S_4(\vp,\epsilon)\le
\frac{2(\lambda-\delta+\epsilon)}{1+\lambda}.
$$
Notice that this bound is smaller than
$\frac{2-2\lambda}{1+\lambda}+3\epsilon$ for
$\delta>2\lambda-1-\epsilon$. As $\epsilon$ can be chosen
arbitrarily small, we derive Theorem~\ref{th1} for $\delta\ge
2\lambda-1$. Since $\delta$ is always nonnegative, the proof is
completed for $\lambda\le \frac12$. Otherwise, we may assume that
$\lambda>\frac12+\frac12 \epsilon$ and the parameter $\delta$ in the
parameter space satisfies $\delta\le 2\lambda-1-\epsilon$. This
concludes the step~8.1 of the proof.

Also notice that for all $x\in J$,~\eqref{eq1} gives
$$
H \asymp |P'''(x)| \ll Q^{-\frac{1+\lambda}{2}+3\eta},
$$
therefore we also have $\eta \gtrsim -\frac{1+\lambda}{6}$.

We apply one more discrete partitioning and split the remaining set
$\DDD_4(k)$ into $\DDD_5(k)\cup \DDD_6(k)\cup \DDD_7(k)$ where for
each $\vq\in \DDD_5(k)$ there exists a rational number $u/v$ with
$|v|\le Q^{\frac{1-\lambda}{3}}$ such that
\begin{equation}\label{eq26}
\left|\frac{q_1}{q_0} - \frac{u}{v}\right| \le
Q^{-\frac{1+\lambda}{3}}.
\end{equation}
Similarly, for each $\vq\in \DDD_6(k)$ there exists $u/v\in\QQ$ with
$|v|\le Q^{\frac{1-\lambda}{4}}$ such that
\begin{equation}\label{eq27}
\left|\frac{q_1}{q_0} - \frac{u}{v}\right| \le
Q^{-\frac{1+\lambda}{4}}.
\end{equation}
Finally, $\DDD_7(k)$ consists of all the remaining points.

For $i\in\{5,6\}$ the Hausdorff dimension of $S_i$ is easy to
compute. The set $\bigcup_{\vq\in \DDD_i(k)} B_\vq$ is covered by
the intervals $B\left(u/v, Q^{-\frac{1+\lambda}{i-2}}\right)$. And
we have at most $\ll Q^{\frac{2(1-\lambda)}{i-2}}$ of them that
intersect $I$. Then
$$
\sum_{k=1}^\infty \#\DDD_i(k) Q^{-\frac{s(1+\lambda)}{i-2}}<\infty
$$
as soon as
$$
\frac{2(1-\lambda)}{i-2}< \frac{s(1+\lambda)}{i-2}\quad
\Longleftrightarrow\quad s\ge \frac{2-2\lambda}{1+\lambda}.
$$
Therefore $\dim S_i < \frac{2-2\lambda}{1+\lambda}$. This concludes
the step~8.3 of the proof.

\medskip
We finish this section with a stronger version of Lemma~4
from~\cite{badziahin_2025} that provides better estimates
than~\eqref{eq13}.
\begin{lemma}\label{lem7}
Suppose that $\Delta_m\cap \ZZ^{4}$ contains at least three linearly
independent vectors. Then
\begin{equation}\label{lem7_eq}
\#(\Delta_m\cap \ZZ^4)\ll Q^{\frac{3-5\lambda}{2}+\delta_m}.
\end{equation}
If $\Delta_m\cap \ZZ^4$ contains at most two linearly independent
vectors then
\begin{equation}\label{lem7_eq2}
\#(\Delta_m\cap \ZZ^4)\ll
Q^{\frac{3-5\lambda}{2}+\delta_m+\delta^*_m},
\end{equation}
where $Q^{\delta^*_m}$ is the value of the third successive minimum
of $\Delta_m$.
\end{lemma}

\proof Let $\tau_1, \ldots, \tau_{n+1}$ be successive minima of
$\Delta_m$. Then from Minkowski's second theorem we get
$$
\prod_{i=1}^{n+1} \tau_i \asymp (\vol(\Delta_m))^{-1} \asymp
Q^{\frac{5\lambda - 3}{2}}.
$$
On the other hand,
$$
\#(\Delta_m\cap \ZZ^4) \asymp \prod_{\tau_i <1} \tau_i^{-1} \asymp
Q^{\frac{3-5\lambda}{2}} \prod_{\tau_i\ge 1} \tau_i.
$$
Then the lemma immediately follows.
\endproof

\section{The proof of steps 8.4 -- 8.6}

From now on we assume that all points $\vq\in \DDD_7(k)$. In view of
the inequality $\eta\ge -\frac{1+\lambda}{6}$ we derive that for all
$\vq$ in this set, the corresponding interval $J$ has length at most
$Q^{-\frac{1+\lambda}{3}}$ and therefore all the rational points
$u/v\in J$ satisfy
\begin{equation}\label{eq25}
|v|\ge Q^{\frac{1-\lambda}{3}}.
\end{equation}

Notice that the interval $J$ must contain one of the roots $b/a$ or
$d/c$ of $P_\va$. Also, if we have several intervals $J(m_1)$,
$J(m_2)$ that contain the same rational point $u/v$ we can cover all
of them by one interval of size
$$
J^*:=\left[\frac{u}{v}-cQ^{\frac{1+\lambda}{2}-\eta},
\frac{u}{v}+cQ^{\frac{1+\lambda}{2}-\eta}\right]
$$
for some absolute constant $c$, and count the number of points
$\vq\in \DDD_7(k)$ lying inside this potentially bigger $J^*$
instead of separately counting these numbers for each $J(m_1)$,
$J(m_2)$, etc. Suppose that the interval $J$ does not contain
rational numbers $u/v$ with $v\ll H^{1/3}$. Since $a^2c \le H$ that
implies that $J$ can only contain one of $b/a$ or $d/c$ but not
both. Hence we can discretely partition $\DDD_7(k)$ into three
subsets $\DDD_8(k), \DDD_9(k)$ and $\DDD_{10}(k)$ where for $\vq\in
\DDD_8(k)$ the corresponding intervals $J$ contain $d/c$ but not
$b/a$. For the set $\DDD_9(k)$ the corresponding intervals $J$
contain $b/a$ but not $d/c$. And finally, for $\DDD_{10}(k)$ the
corresponding intervals $J$ contain rational points $u/v$ with $v\ll
H^{1/3}$.

{\bf The case $d/c\in J$, $b/a\not\in J$.} By examining the
derivative $P'_\va(x)$ we find that the largest value of
$|P_\va(x)|$ for $x$ between $b/a$ and $d/c$ is for $x_0 =
\frac{b}{3a} + \frac{2d}{3c}$ and one can quickly check that
$|x_0-b/a|\asymp |d/c-b/a|$. For $x_0$ we must have $|P_\va(x_0)| >
HQ^{-1-\lambda}$ therefore there exists $x\in J$ such that
$$
\left|x-\frac{d}{c}\right|\gg Q^{-\frac{1+\lambda}{2} -
\eta},\quad\mbox{and}\quad \left|x-\frac{b}{a}\right| \asymp
\left|\frac{d}{c}-\frac{b}{a}\right|.
$$
Since $|P_\va(x)| = H|x-d/c||x-b/a|^2 < HQ^{-1-\lambda}$, the above
inequalities imply
\begin{equation}\label{eq16}
\left|\frac{d}{c}- \frac{b}{a}\right| \ll
Q^{-\frac{1+\lambda-2\eta}{4}}.
\end{equation}
The left hand side is always at least $\frac{1}{ac} \ge H^{-1}$.
Then for $\eta<0$, in view of~\eqref{eq12} we derive
$$
Q^{-\lambda+\delta-\epsilon}\ll \frac{1}{H}\ll
Q^{-\frac{1+\lambda-2\eta}{4}}
$$
or
\begin{equation}\label{eq17}
\eta \gtrsim \frac{1-3\lambda}{2} + 2\delta - 2\epsilon.
\end{equation}

Recall that we have $|c|>H^{1/3}$, therefore $a\le (H/c)^{1/2} \le
H^{1/3}$. Now, for a fixed $a$ the value of $c$ can not exceed
$|c|\le H/a^2$. The number of fractions $d/c$ that lie in an
interval of length $l$ with $1\le c\le C$ can be estimated as $\ll
\max\{1,C^2l\}$. Therefore, for a fixed $b/a$, the number of
fractions $d/c$ that satisfy~\eqref{eq16}, is bounded from above by
$\ll\max\{1,H^2a^{-4}Q^{-\frac{1+\lambda-2\eta}{4}}\}$. Summing over
all fractions with denominator $a$ and then over all $a$, we end up
with the following bound on the number of intervals $J$ with $d/c\in
J$ and $b/a\not\in J$:
$$
\sum_{a=1}^{H^{1/3}} \sum_b
\max\left\{1,\frac{H^2}{a^4}Q^{-\frac{1+\lambda-2\eta}{4}}\right\}
\ll H^{2/3}+ H^2 Q^{-\frac{1+\lambda-2\eta}{4}}.
$$
With help of~\eqref{eq12} and~\eqref{eq13} we get that the total
number of the corresponding points $\vq\in\DDD_8(k,\vp,\epsilon)$
for $\eta\ge 0$ is bounded from above by
$$
\ll Q^{\frac23(\lambda-\delta-\eta+\epsilon) +
\frac{3-5\lambda}{2}+2\delta+2\epsilon} + Q^{\frac{7\lambda-1}{4} -
\frac32 \eta-2\delta +2\epsilon + \frac{3-5\lambda}{2} +
2\delta+2\epsilon}.
$$
This sum maximises when $\eta$ is the smallest possible, i.e.
$\eta=0$, and $\delta$ is the largest possible, i.e. $\delta =
2\lambda-1$. Then the estimate becomes
$$
\#\DDD_8(k,\vp,\epsilon)\ll Q^{\frac{5\lambda+1}{6} +
\frac83\epsilon} + Q^{\frac{5-3\lambda}{4} + 4\epsilon}
\stackrel{\scriptscriptstyle \lambda\le 3/5}\ll Q^{2-2\lambda +
4\epsilon}.
$$

For $\eta<0$, the inequality~\eqref{eq16} together with the fact
that all $\vq$ do not lie in $\DDD_6(k)$ implies that $a\gg
Q^{\frac{1-\lambda}{4}}$. That modifies the number of possible
intervals $J$ to
$$
\sum_{a=c_4Q^{\frac{1-\lambda}{4}}}^{H^{1/3}} \sum_b
\max\left\{1,\frac{H^2}{a^4}Q^{-\frac{1+\lambda-2\eta}{4}}\right\}
\ll H^{2/3}+ H^2 Q^{-\frac{3-\lambda-2\eta}{4}}
$$
for some absolute constant $c_4>0$. Then~\eqref{eq12}
and~\eqref{eq13} give the upper bound for the number of points
$\vq\in \DDD_8(k,\vp,\epsilon)$ as
\begin{equation}\label{eq18}
\# \DDD_8(k,\vp,\epsilon) \ll Q^{\frac23(\lambda-\delta+\epsilon) +
\frac{3-5\lambda}{2}+2\delta-\eta+3\epsilon} + Q^{\frac{9\lambda -
3}{4} -2\delta+\frac{\eta}{2} + 2\epsilon +\frac{3-5\lambda}{2} +
2\delta-\eta+3\epsilon}
\end{equation}
The first summand maximises when $\eta$ is as small as possible,
which in view of ~\eqref{eq17} gives that it is at most
$$
Q^{\frac{9-11\lambda}{6} +\frac43 \delta + \frac{11}{3}\epsilon
-\frac{1-3\lambda}{2}-2\delta + 2\epsilon}.
$$
Now this expression maximises when $\delta=0$ and then we get
$$
\ll Q^{\frac{3-\lambda}{3} + \frac{14}{3}\epsilon}
\stackrel{\lambda\le 3/5}\ll Q^{2-2\lambda + 6\epsilon}.
$$

The second summand in~\eqref{eq18} maximises when $\eta$ is smallest
possible and it does not depend on $\delta$. Here we can use a
weaker estimate $\eta\ge -\frac{1+\lambda}{6}$. By substituting this
into the summand we get that it is at most
$$
\ll Q^{\frac{5-\lambda}{6} + 5\epsilon} \stackrel{\lambda\le 3/5}\ll
Q^{2-2\lambda + 5\epsilon}.
$$

We conclude that in all the cases we have
$$
\#\DDD_8(k,\vp,\epsilon) \ll Q^{2-2\lambda + 6\epsilon}.
$$

\medskip
\medskip{\bf The case $b/a\in J, d/c\not\in J$.} Notice that for
$\big|x-\frac{b}{a}\big|\le Q^{-\frac{1+\lambda}{2}}$ we have
$|P_\va(x)|\le HQ^{-1-\lambda}$. Therefore in this case we must have
$\eta\le 0$. On the other hand, the bound~\eqref{eq1} for $x_0 =
\frac{b}{a}$ gives
\begin{equation}\label{eq19}
\frac{1}{ac}\le \left|\frac{d}{c}-\frac{b}{a}\right| \asymp
\frac{|P''_\va(x_0)|}{H(P_\va)}\ll Q^{2\eta}.
\end{equation}

We continuously partition the set $\DDD_9(k)$ by adding the
parameter $\sigma$ such that $|c| = Q^\sigma$ to the family of
$\vp$. Obviously, $\sigma\ge 0$. On the other hand, since $c\le H$,
we have $\sigma\lesssim \lambda-\delta+\epsilon$. Now we compute the
number of intervals $J$ in two different ways. Firstly, since we
have $|a|\le (H/c)^{1/2} \le (H/Q^{\sigma})^{1/2}$, the number of
different fractions $b/a\in I$ (and hence the number of intervals
$J$) is bounded from above by
$$
\#J \ll H/Q^{\sigma}\ll Q^{\lambda - \delta-\sigma + \epsilon}.
$$

Secondly, for a given fraction $d/c$, the number of fractions $b/a$
with $a\le (H/c)^{1/2}$ that satisfy~\eqref{eq19}, is bounded from
above by $\max\{1,HQ^{2\eta}/c\}$. But we also have $Q^{2\eta}\ge
(ac)^{-1}$ which gives $HQ^{2\eta}/c\gg \frac{a}{c}$, so for $a\ge
H^{1/3}$ the second term in the maximum prevails. Summing over all
possible fractions $d/c$ that have $c\le Q^{\sigma+\epsilon}$ gives
$$
\# J\ll \sum_{c\ge Q^{\sigma}}^{Q^{\sigma+\epsilon}} \sum_d
\max\left\{1,\frac{H}{c}Q^{2\eta}\right\} \ll HQ^{\sigma +
2\eta+\epsilon}.
$$

Suppose that $\eta\ge \delta-\frac{1-\lambda}{2}-\sigma$. Then we
use the first estimate for $\#J$ and in view of~\eqref{eq13}, we
have that the number of points $\vq\in \DDD_9(k,\vp,\epsilon)$ is
$$
\#\DDD_9(k,\vp,\epsilon)\ll Q^{\lambda - \delta-\sigma + \epsilon +
\frac{3-5\lambda}{2} + 2\delta - \eta +3\epsilon} =
Q^{\frac{3}{2}(1-\lambda) + \delta - \eta - \sigma + 4\epsilon} \ll
Q^{2-2\lambda + 4\epsilon}.
$$
Next, suppose that $\eta\le \delta-\frac{1-\lambda}{2}-\sigma$. Then
we apply the second estimate for $\#J$. Again, apply~\eqref{eq13} to
get the upper bound for the number $\vq\in \DDD_9(k,\vp,\epsilon)$:
$$
\#\DDD_9(k,\vp,\epsilon)\ll Q^{\lambda - \delta + \epsilon
+\sigma+2\eta + \epsilon+\frac{3-5\lambda}{2} + 2\delta -\eta +
3\epsilon}.
$$
This estimate maximises when $\eta$ and $\delta$ are maximal
possible, i.e. $\delta = 2\lambda - 1$ and $\eta = \delta -
\frac{1-\lambda}{2}-\sigma$. That gives
$$
\ll Q^{\frac{3-3\lambda}{2} -\frac{1-\lambda}{2} +2\delta
+5\epsilon} \ll Q^{3\lambda - 1 + 5\epsilon} \stackrel{\lambda\le
3/5}\ll Q^{2-2\lambda + 5\epsilon}.
$$
In all cases we get $\#\DDD_9(k,\vp,\epsilon)\ll Q^{2-2\lambda +
5\epsilon}$.

To finish the steps~8.5 and~8.6 of the proof we notice that for both
sets $\DDD_8(k,\vp,\epsilon)$ and $\DDD_9(k,\vp,\epsilon)$ we can
take the notion $N(k,\vp,\epsilon)$ from Lemma~\ref{lem3} to be
$Q^{2-2\lambda+6\epsilon}$ and then
$$
M(k,\epsilon) = \frac{2-2\lambda + 6\epsilon}{1+\lambda}\;
\stackrel{\epsilon\to 0}\to\; \frac{2-2\lambda}{1+\lambda}.
$$

\section{Fractions with small height inside J: the proof of steps 8.7.1 -- 8.7.4}

We now focus on the set $\DDD_{10}(k)$. Recall that for all $\vq$ in
this set the corresponding interval $J(\vq)$ contains a rational
number $u/v$ with $v\le H^{1/3}$. Then the number of such fractions
$u/v\in I$ and hence the number of intervals $J(\vq)$ is bounded
from above by $H^{2/3}$. Also notice that since $\eta\ge
-\frac{1+\lambda}{6}$, the length of $J(\vq)$ satisfies $|J| \ll
Q^{-\frac{1+\lambda}{3}}$ and since $\vq\not\in \DDD_5(k)$, we must
have
\begin{equation}\label{eq33}
H^{1/3}\gg v \gg Q^{\frac{1-\lambda}{3}}\quad \Longrightarrow\quad H
\gg Q^{1-\lambda}.
\end{equation}

We continuously partition $\DDD_{10}(k)$ by adding a new parameter
$\delta^*$ to the list of parameters $\vp$. It is defined by the
third successive minimum $\tau_3 = Q^{\delta^*}$ of the box
$\Delta_m$. Since $\delta^*\le \delta$, it is definitely bounded
from above. One can use the second Minkowski theorem about the
successive minima of $\Delta_m$ to see that $\delta^*$ is bounded
from below as well.

Consider the case when $\delta^* + \delta\le
\frac12(1-\lambda)+\epsilon$. Then, with help of Lemma~\ref{lem7}
the number of points $\vq\in \DDD_{10}(k,\vp, \epsilon)$ that
correspond to those parameters $\delta, \delta^*$ and lie in a given
interval $J$, is bounded from above by
$$
\ll Q^{\frac{3-5\lambda}{2} + \delta + \delta^* - \eta+2\epsilon}.
$$
Summing up over all intervals $J$, we estimate
$$
\#\DDD_{10}(k,\vp,\epsilon)\ll Q^{\frac23(\lambda-\delta+\epsilon) +
\frac{3-5\lambda}{2} + \delta + \delta^* - \eta+2\epsilon} =
Q^{\frac{9-11\lambda}{6} + \frac13\delta+\delta^* - \eta +
\frac83\epsilon}.
$$
If $\delta+\delta^*$ is fixed, this expression maximises when
$\delta^*$ is maximal possible, i.e. $\delta^* = \delta$. We also
need to take $\eta$ as small as possible to maximise the expression
\big(i.e. $\eta = -\frac{1+\lambda}{6}$\big). That gives us
\begin{equation}\label{eq30}
\#\DDD_{10}(k,\vp,\epsilon)\ll Q^{\frac{9-11\lambda}{6}
+\frac{1-\lambda}{3} + \frac43\epsilon +  \frac{1+\lambda}{6} +
3\epsilon} \ll Q^{2-2\lambda + 5\epsilon}.
\end{equation}

For the rest of the paper we will assume that
\begin{equation}\label{eq32}
\delta+\delta^*\ge \frac12(1-\lambda)+\epsilon.
\end{equation}
Notice that this inequality also implies that $\delta^*\ge 0$
because otherwise we have $\delta \ge \frac12(1-\lambda)+\epsilon >
2\lambda-1$ which is a contradiction with $\lambda\le 3/5$. In
particular, this means that $\Delta_m \cap
\DDD_{10}(k,\vp,\epsilon)$ lies in a two-dimensional subspace.

%Similar computations show that for $\delta\ge \frac{1-\lambda}{4}$,
%$$
%\#\vq \ll Q^{\frac23(\lambda - \delta+ \epsilon) +
%\frac{3-5\lambda}{2} + 2\delta - \eta + 2\epsilon} \ll
%Q^{\frac{10-10\lambda}{6} + \frac43\cdot \frac{1-\lambda}{4} +
%3\epsilon} \ll Q^{2-2\lambda + 3\epsilon}.
%$$

Consider the boxes $\Delta_m$ such that for any two distinct
primitive points $\vq_1,\vq_2\in \Delta_m\cap
\DDD_{10}(k,\vp,\epsilon)$ one has
\begin{equation}\label{eq28}
\left|\frac{q_{11}}{q_{10}} - \frac{q_{21}}{q_{20}}\right| \ge
Q^{-\frac{4-2\lambda}{3} - \epsilon}.
\end{equation}
In this case, the number of points $\vq\in
\DDD_{10}(k,\vp,\epsilon)\cap \Delta_m$ is bounded from above by
$$
Q^{-\frac{1+\lambda}{2}+\frac{4-2\lambda}3 + \epsilon}\le Q^{\frac{5
- 7\lambda}{6}+\epsilon}.
$$
Discretely partition $\DDD_{10}(k)$ into $\DDD_{11}(k)$ and
$\DDD_{12}(k)$ where all $\vq\in\DDD_{11}(k)$ lie in one of the
boxes $\Delta_m$ that satisfy the above condition, and
$\DDD_{12}(k)$ contains the remaining points. Summing over
$\Delta_m$ that comprise $\DDD_{11}(k)$ and correspond to a given
interval $J$ and then summing over all intervals $J$, the following
upper bound for $\#\DDD_{11}(k,\vp,\epsilon)$ is satisfied:
$$
\#\DDD_{11}(k,\vp,\epsilon) \ll Q^{\frac23(\lambda - \delta +
\epsilon)+\frac{5-7\lambda}{6}+\max\{0,-\eta\} + \epsilon}.
$$
This sum maximises when both $\delta$ and $\eta$ are smallest
possible, i.e. $\delta=0$ and $\eta = -\frac{1+\lambda}{6}$. Then we
have
\begin{equation}\label{eq31}
\#\DDD_{11}(k,\vp,\epsilon)\ll Q^{\frac{6-2\lambda}{6} + 2\epsilon}
\stackrel{\lambda\le 3/5}\ll Q^{2-2\lambda + 2\epsilon}.
\end{equation}

For the remaining boxes $\Delta_m$ we have
$\DDD_{12}(k,\vp,\epsilon)\cap \Delta_m$ lie in a two-dimensional
subspace with $\delta_m^*+\delta_m\ge \frac{1-\lambda}{2}+\epsilon$
and there exist two points $\vq_1,\vq_2\in
\DDD_{12}(k,\lambda,\epsilon)\cap \Delta_m$ with
$$|q_{11}/q_{10} -
q_{21}/q_{20}|< Q^{-\frac{4-2\lambda}{3}-\epsilon}.$$ Notice that in
view of the second Minkowski second theorem, the successive minima
of $\Delta_m$ satisfy
$$
\prod_{i=1}^4 \tau_i \asymp Q^{\frac{5\lambda-3}{2}}\le 1,
$$
therefore we have
\begin{equation}\label{eq22}
\tau_1\tau_2\ll Q^{-\frac{1-\lambda}{2}-\epsilon}.
\end{equation}

Consider two points
$$
\vp_1 := (v^3, v^2u, vu^2, u^3),\quad \vp_2 := \left(v^2r, v(ur+1),
u(ur+2), \frac{u^2(ur+3)}{v}\right),
$$
where $0\le r<v$ is such that $v\mid ur+3$. One can easily check
that both points satisfy the equations
$$
-x_0^2 p_0 + 2x_0 p_1 - p_2 = -2x_0^2p_0 + 3x_0^2p_1 - p_3 = 0
$$
where $x_0 = \frac{u}{v}$, therefore $\Span(\vp_1, \vp_2)$ always
lies inside the hyperplane $\va\cdot\vx=0$ with the corresponding
polynomial $P_\va(x) = (vx-u)^2$. Finally, these two points are
obviously linear independent.

Consider a point $x = x_0 + t\in J$ which is the center of one of
the boxes $\Delta_m\subset J$. Because of~\eqref{eq14} and the fact
that $\eta\ge -\frac{1+\lambda}{6}$, we have $|t|\le
Q^{-\frac{1+\lambda}{3}}$. We compute
\begin{equation}\label{eq23}
|p_{10} x- p_{11}| = |v^3(x_0+t)-v^2u| = |t v^3|\ll
HQ^{-\frac{1+\lambda}{3}};
\end{equation}
\begin{equation}\label{eq24}
|p_{20}x - p_{21}|=  |v^2r(x_0 + t) - v(ur+1)| \le |v| + |tv^2r|\ll
H^{1/3} + HQ^{-\frac{1+\lambda}{3}}.
\end{equation}
We also compute
%For $\lambda\le \frac35$ we get
%$$
%H^{1/3}\ll Q^{\lambda/3} \le Q^{\frac{1-\lambda}{2}}
%$$
%and
%$$
%HQ^{-\frac{1+\lambda}{3}} \ll Q^{\frac{2\lambda-1}{3}} \le
%Q^{\frac{1-\lambda}{2}}.
%$$
%Therefore, for both points $\vq=\vq_0$ and $\vq=\vq_1$ we have
%$|q_0x-q_1|\ll Q^{\frac{1-\lambda}{2}}$. Next, compute
\begin{equation}\label{eq20}
|-v^3 (x_0+t)^2 + 2v^2u (x_0+t) - vu^2| = |v^3t^2|\ll
HQ^{-\frac{2(1+\lambda)}{3}};
\end{equation}
\begin{equation}\label{eq21}
|-v^2r(x_0+t)^2 + 2v(ur+1)(x_0+t) - u(ur+2)| = |2vt - v^2rt^2| \le
H^{1/3}Q^{-\frac{1+\lambda}{3}} + HQ^{-\frac{2(1+\lambda)}{3}}.
\end{equation}
Notice that $HQ^{-\frac{2(1+\lambda)}{3}} < H^{1/3}
Q^{-\frac{1+\lambda}{3}}$ is equivalent to $H<
Q^{\frac{1+\lambda}{2}}$ which is true due to~\eqref{eq12} and
$\lambda<1$. Finally, by analogous computations one derives the same
bound~\eqref{eq20} for $|-2x^3 p_{10} + 3x^2 p_{11} - p_{13}|$ and
the bound~\eqref{eq21} for $|-2x^3 p_{20} + 3x^2 p_{21} - p_{23}|$.

Choose $\vq_1,\vq_2\in \DDD_{12}(k,\vp,\epsilon)\cap \Delta_m$ such
that $|q_{11}/q_{10}-q_{21}/q_{20}| < Q^{-\frac{4-2\lambda}{3} -
\epsilon}$. Suppose that $\vp_1,\vq_1$ and $\vq_2$ are linearly
independent and consider the hyperplane $\PPP$ that passes through
these points. Since $\DDD_{12}(k,\vp,\epsilon)\cap \Delta_m$ lies in
a two-dimensional space and $\vq_1$, $\vq_2$ are linearly
independent, $\PPP$ goes through all points from
$\DDD_{12}(k,\vp,\epsilon)\cap \Delta_m$ and hence its height should
be at least $H$. On the other hand, it is bounded by $|\vp_1\wedge
\vq_1\wedge\vq_2|$. We compute
$$
\begin{vmatrix}
p_{10}&p_{11}&p_{12}\\
q_{10}&q_{11}&q_{12}\\
q_{20}&q_{21}&q_{22}\\
\end{vmatrix} = \begin{vmatrix}
p_{10}&p_{11}-p_{10}x&p_{12}-2p_{11}x+p_{10}x^2\\
q_{10}&q_{11}-q_{10}x&q_{12}-2q_{10}x+q_{10}x^2\\
q_{20}&q_{21}-q_{20}x&q_{22}-2q_{20}x+q_{20}x^2\\
\end{vmatrix}.
$$
For $x = q_{11}/q_{10}$ we derive the following upper bounds on the
absolute values of each entry of this matrix:
$$
\begin{vmatrix}
H& H|t| & H|t|^2\\
Q& 0& Q^{-\lambda}\\
Q& Q^{\frac{2\lambda-1}{3}-\epsilon}&Q^{-\lambda}
\end{vmatrix}
$$
By examining all the terms in the determinant, we derive that it is
bounded from above by $$H(Q^{1-\lambda}|t| +
Q^{\frac{2\lambda+2}{3}-\epsilon}|t|^2+
Q^{-\frac{\lambda+1}{3}-\epsilon}).$$ in view of $|t|\ll
Q^{-\frac{1+\lambda}{3}}$, we get that the first term is less than
$Q^{-\frac23\epsilon}$ as soon as $\lambda \ge
\frac12+\frac12\epsilon$. Recall that for $\lambda\le
\frac12+\frac12\epsilon$ we have already establish the theorem, so
without loss of generality we can assume this condition. The second
summand is at most $Q^{-\epsilon}$. We conclude that this
determinant, which gives one of the coordinates of $|\vp_1\wedge
\vq_1\wedge\vq_2|$, is at most $HQ^{-\epsilon}$. Analogous
computations give that the determinants of the following coordinates
of this multivector are also less than $HQ^{-\epsilon}$;
$$
\begin{vmatrix}
p_{10}&p_{11}&p_{13}\\
q_{10}&q_{11}&q_{13}\\
q_{20}&q_{21}&q_{23}\\
\end{vmatrix}; \qquad \begin{vmatrix}
p_{10}&p_{12}&p_{13}\\
q_{10}&q_{12}&q_{13}\\
q_{20}&q_{22}&q_{23}\\
\end{vmatrix}.
$$
Finally, the remaining coordinate of $$|\vp_1\wedge
\vq_1\wedge\vq_2|$$ also has to be smaller than $\ll
HQ^{-\epsilon}$. This contradicts the fact that $||\vp_1\wedge
\vq_1\wedge\vq_2||\ge H$. We conclude that $\vp_1\in
\Span(\vq_1,\vq_2) = \Span(\Delta_m\cap \DDD_{12}(k,\vp,\epsilon))$.

Consider the points $\vq\in \Delta_m$ that realises the first
successive minimum of $\Delta_m$, i.e.
$$
|q_0|\le \tau_1 Q,\; |q_0x_m - q_1|\le \tau_1
Q^{\frac{1-\lambda}{2}},\; |(1-i)x_m^iq_0 + ix_m^{i-1} - q_i|\le
\tau_1 Q^{-\lambda},\; i\in\{2,3\}.
$$
Suppose that $\vp_1,\vp_2$ and $\vq$ are linearly independent and
estimate $||\vp_1\wedge \vp_2\wedge\vq||$. Proceeding as before and
using~\eqref{eq23},~\eqref{eq24},~\eqref{eq20} and~\eqref{eq21}, we
compute the following upper bounds for the entries of the following
determinant
$$
\begin{vmatrix}
p_{10}&p_{11}&p_{12}\\
p_{20}&p_{21}&p_{22}\\
q_{0}&q_{1}&q_{2}\\
\end{vmatrix} \ll \begin{vmatrix}
H&H|t|&H|t^2|\\
H&H^{1/3}&H^{1/3}|t|\\
\tau_1 Q& \tau_1Q^{\frac{1-\lambda}{2}}&\tau_1Q^{-\lambda}
\end{vmatrix}
$$
By expanding this determinant and estimating all the terms, we get
the upper bound
\begin{equation}\label{eq35}
\tau_1H(QH^{1/3}|t|^2 + Q^{\frac{1-\lambda}{2}}H^{1/3}|t| +
Q^{\frac{1-\lambda}{2}}H|t|^2).
\end{equation}
Notice that in view of $|t|\ge Q^{-\frac{1+\lambda}{2}}$ and $H\ll
Q^{\frac{1+\lambda}{2}}$, the first term in this sum is bigger than
the others and then we use~\eqref{eq32} and~\eqref{eq22} to
continue estimating the term:
$$
\ll H\cdot Q^{-\frac{1-\lambda}{4}-\frac12\epsilon + 1+
\frac13\left(\lambda-\frac{1-\lambda}{4}+\frac{\epsilon}{2}\right) -
\frac{2(1+\lambda)}{3}} = HQ^{-\frac\epsilon{3}}.
$$
Again, analogous computations for the other two determinants in
$\va=||\vp_1\wedge\vp_2\wedge\vq||$ also give that they are smaller
than $HQ^{-\epsilon/3}$ and the equation $\va\cdot\vq = 0$ gives the
same estimate for the last term in $\va$. We get a contradiction
with $||\vp_1\wedge\vp_2\wedge\vq||\ge H$.

We conclude that $\vp_1,\vp_2,\vq$ are linearly dependent, i.e.
$\vq\in \Span(\vp_1,\vp_2)$. If $\vp_1$ and $\vq$ are linearly
independent then $\DDD_{12}(k,\vp,\lambda)\cap \Delta_m \subset
\Span(\vq,\vp_1) = \Span(\vp_1,\vp_2)$. But both $\vp_1,\vp_2$ lie
in the plane with the corresponding polynomial $P_\va(x) =
(ax-b)^2$, i.e. its height is $H^{2/3}<H$ --- a contradiction.

We are left with the case when $\vq$ is collinear with $\vp_1$, and
since $\vp_1$ is primitive, it realises the successive minimum
$\tau_1$. Therefore we have
$$
\tau_1 \asymp \max\left\{\frac{H}{Q},\;
\frac{H|t|}{Q^{\frac{1-\lambda}{2}}},\; H|t|^2Q^\lambda\right\}.
$$
In view of $|t|\ge Q^{-\frac{1+\lambda}{2}}$, we observe that
$\tau_1\asymp H|t|^2Q^\lambda$.

Next, consider the point $\vq\in \Delta_m$ that realises the second
successive minimum of $\Delta_m$. It is now linearly independent
with $\vp_1$ by construction. If $\vp_1,\vp_2,\vq$ are linearly
dependent then $\DDD_{12}(k,\vp,\lambda)\cap \Delta_m \subset
\Span(\vq,\vp_1) = \Span(\vp_1,\vp_2)$ which leads to a
contradiction as both $\vp_1$ and $\vp_2$ lie in the hyperplane of
height $H^{2/3}<H$. Therefore we get that
$||\vp_1\wedge\vp_2\wedge\vq||\gg H$. This case is considered
analogously to the previous one and we derive that the length of
this multivector is bounded from above by~\eqref{eq35} with $\tau_2$
in place of $\tau_1$. Then we must have
$$
\tau_2 QH^{1/3}|t|^2 \gg 1.
$$
Notice that $H|t|^2 Q^\lambda \asymp \tau_1$ and therefore we get
$$
\tau_1\tau_2 Q^{1-\lambda} H^{-2/3} \gg 1\quad
\Longleftrightarrow\quad \tau_1\tau_2 \gg H^{2/3}Q^{-(1-\lambda)}.
$$
On the other hand, by~\eqref{eq22} we have $\tau_1\tau_2 \ll
Q^{-\delta-\delta^*}\ll Q^{-\frac{1-\lambda}{2}}$. By combining this
upper bound with the lower bound above, we derive
$$
H\ll Q^{\frac34 (1-\lambda)}.
$$
But this contradicts~\eqref{eq33}. We finally exhaust all the cases
and therefore the set $\DDD_{12}(k)$ is empty. That establish the
step~8.7.4.

To accomplish steps 8.7.1 -- 8.7.3, we notice that for both
$\DDD_{10}(k,\vp,\epsilon)$ and $\DDD_{11}(k,\vp,\epsilon)$ we can
take the notion $N(k,\vp,\epsilon)$ from Lemma~\ref{lem3} to be
$Q^{2-2\lambda+5\epsilon}$ (see~\eqref{eq30} and~\eqref{eq31}) and
hence
$$
M(k,\epsilon) = \frac{2-2\lambda+5\epsilon}{1+\lambda}\;
\stackrel{\epsilon\to 0}\to\; \frac{2-2\lambda}{1+\lambda}.
$$
The proof of Theorem~\ref{th1} is now complete.

%The upshot is, in all the cases when we do not get a contradiction,
%we derive
%$$
%\# \DDD_2(k,\lambda,\epsilon) \ll Q^{2-2\lambda + 5\epsilon}.
%$$
%Then we can choose $N(k,\vp,\epsilon)$ (see Lemma~\ref{lem3}) such
%that $\sup_{\vp\in \FFF} N(k,\vp,\epsilon) = Q^{2-2\lambda +
%5\epsilon}$. then
%$$
%M(k,\epsilon) = \frac{\log(\sup_{\vp\in \FFF}
%N(k,\vp,\epsilon))}{\log \rho^+(k)} =
%\frac{2-2\lambda+5\epsilon}{1+\lambda}.
%$$
%This expression tends to $\frac{2-2\lambda}{1+\lambda}$ uniformly as
%$\epsilon\to 0$, The proof of Theorem is now complete.

\bigskip
\noindent Dzmitry Badziahin\\ \noindent The University of Sydney\\
\noindent Camperdown 2006, NSW (Australia)\\
\noindent {\tt dzmitry.badziahin@sydney.edu.au}

\end{document}